\documentclass{article}

\begin{document}

\section*{Model problems for two equations, which type depends on solution}


\begin{center}
Isaac I. Vainshtein \\ isvain@mail.ru\\
Department of applied mathematics and computer security \\
The institute of Space and information technologies,\\
Siberian Federal University,\\
Russian Federation
\end{center}

\begin{abstract}
In this work there are considered model problems for
two nonlinear equations, which type depends on the solution. One of
the equations may be called a nonlinear analog of the
Lavrent'ev-Bitsadze equation.
\end{abstract}

Keywords: type of equation, equations of elliptic, hyperbolic type,
the Lavrent'ev-Bitsadze equation, Tricomi problem.

\section*{Introduction }


In this work there are considered two nonlinear equations that may be models for research and formulation of boundary value problems for nonlinear equations, which type depends on the solution:
\begin{equation}
\frac{\partial^2U(x,y)}{\partial x^2}+\textrm{sign}\,
  U \frac{\partial^2U(x,y)}{\partial y^2}=0,
\end{equation}
\begin{equation}
\left(\frac{\partial^2U(x,y)}{\partial
x^2}\right)^2-\left(\frac{\partial^2U(x,y)}{\partial
y^2}\right)^2= \left(\frac{\partial^2U}{\partial
x^2}+\frac{\partial^2U}{\partial
y^2}\right)\left(\frac{\partial^2U}{\partial
x^2}-\frac{\partial^2U}{\partial y^2}\right)=0.
\end{equation}

The equation (1) changes type, if the solution switches sign or even if vanishes in considered domain.
When $U>0$, we have Laplace equation (elliptic type). When $U<0$, we have wave equation (hyperbolic type).
If solution $U(x,y)$ vanishes in $B\subset D$ set, there is parabolic degeneracy in set $B$. $D$ is the domain of the considered specific solution $U(x,y)$.

The equation (1) may be considered as a nonlinear analogue of Lavrent'ev–Bitsadze equation ~[1]
\begin{equation}
\frac{\partial^2U}{\partial x^2}+\textrm{sign} \,
  y \frac{\partial^2U}{\partial y^2}=0,
\end{equation}
which is  one of the canonical equations in the mixed type equations theory.

The equation type (2) depends on factor that makes solution vanished.

Let us consider the matter of formulation of the boundary value problem for equations (1),(2).
Let $D$ be a domain in which solution is sought.
It should be noted, that solving might be implemented along with the domain search.
In what follows the solution is sought in $C^1(D)\cap
C^2(D_i),\ D_i\subseteq D.$ In domain $D_i$ solution saves equation type: elliptic or hyperbolic.

Solution to the problem that does not change the type of equation in the whole domain $D$
will be called \textit{trivial}. If in the entire area $D$ the solution satisfies the wave equation,
it is the trivial hyperbolic solution. If in the entire area $D$
it satisfies the Laplace equation -- \textit{trivial elliptical one}. Trivial solutions of equation (1)
don't vanish in domain $D$. The solution that satisfies Laplace equation in one part of the domain
and wave equation in another will be called \textit{a changing the type of the equation}.

When setting boundary conditions for the equation  (2), it should be taken into account
the possible or prescribed sign of the solution at the boundary, since there
will be determined the sign of solutions by the continuity in the boundary vicinity and the type of equation as well,
hence the form of possible boundary conditions, typical for given type.

In this connection, the components which define the boundary problem (the area where the solution is sought,
the boundary conditions, the changeable type of equation, smoothness) must be strictly consistent.

In the event that boundary problem is formulated, for example, if the domains with boundary conditions are specified, then it can be considered a problem with unknown line, crossing which the solution changes equation type.

A different approach can exist.  Specify a curve and a domain, which encloses this curve,
that when curve is crossed, the solution will change the equation type and satisfy given boundary conditions.

\section*{1. Model problems.}
Let us consider the number of model problems for (1), (2).

\textbf{Dirichlet problem.} In a bounded domain $D$ with boundary
$\Gamma$ it is required to find the solution of equation satisfying
\begin{equation}
U|_\Gamma=\varphi(s).
\end{equation}
\textbf{Dirichlet problem for equation (1).}

In the case of $\varphi(s)\geq 0$ the harmonic function satisfying the boundary condition (4), is a trivial elliptic solution to the problem, since by the maximum principle the function is positive in $D$. There are no trivial solutions, if boundary function $\varphi(s)$ changes sign.\\
\textbf{Dirichlet problem for equation (2). } A harmonic function that satisfied boundary condition (5)
is a trivial elliptic solution.

\textbf{The model Dirichlet problem.} $D$ is a circle of radius $R$ ($r^2=x^2+y^2\leq R^2$), $\varphi(s)=H.$\\
\textbf{The model Dirichlet problem for equation (1).}

Let $\varphi(s)=H>0.$ The solution $U=H>0$ is trivial elliptic. The function
\begin{equation}
U(x,y)=\left\{ \begin{array}{ll} \displaystyle
\frac{H}{2\ln\frac{a}{R}}(1-\frac{r^2}{a^2}),\quad\mbox{if}\ 0\leq r\leq a, \\ \\
\displaystyle
H\frac{\ln\frac{a}{r}}{\ln\frac{a}{R}},\quad\mbox{if}\ a\leq
r\leq R
\end{array}\right.
\end{equation}
for any fixed $a$ ($0<a<R$) is continuously differentiable in the circle
 $r<R$, negative inside the circle $r<a$ satisfying the wave equation , in the ring $a<r<R$ positive and and harmonic.
At $r=R$, $U=H$, at $r=a$, $U=0$.

Let $\varphi(s)=H<0.$ The function
$$
U(x,y)=C(r^2-R^2)+H
$$
is a trivial hyperbolic solution of Dirichlet problem, $C>\frac{H}{R^2}$ is an arbitrary constant.

Thus, the model Dirichlet problem for equation (1) at $H>0$ has one trivial elliptic along with infinite number of solution that changes it's type. There is infinite number of trivial hyperbolic solutions, when $H<0$.

\textbf{The model Dirichlet problem for equation (2).} There is an infinite number of solutions to this problem:
 \begin{equation}                                                                                                 U(x,y)=\left\{ \begin{array}{ll} \displaystyle                                                                         C(r^2-a^2)+2Ca^2\ln\frac{a}{R}+H,\quad\mbox{if}\ 0\leq r\leq a, \\ \\
\displaystyle 2Ca^2\ln\frac{r}{R}+H,\quad\mbox{if}\ a\leq r\leq
R,
\end{array}\right.
\end{equation}
 that change equation type and infinite number of hyperbolic trivial solutions.
$$
U(x,y)=C(r^2-R^2)+H,
$$
for any $C,a,\ a\in(0,R).$

\textbf{Neumann problem}. In a bounded domain $D$ with boundary
$\Gamma$ it is required to find solution satisfying
\begin{equation}
\left.\frac{\partial U}{\partial n}\right|_\Gamma=\varphi(s).
\end{equation}
Let $\int\limits_\Gamma\varphi(s)ds=0$ è $V(x,y)$ is an arbitrary fixed harmonic function in $D$
that satisfies the boundary condition (7).

There is an infinite number of trivial elliptic solutions to the Neumann problem for equations (1), (2):
$$U=V(x,y)+C,\quad C> \max_D|V(x,y)|$$
for equation (1),
$$U=V(x,y)+C,$$
for equation (2), C -- arbitrary cnostant.

\textbf{The model Neumann problem.}  $D$ --  is a circle of radius $R$, $\varphi(s)=K.$\\
\textbf{The model Neumann problem for equation (1).}

When $K>0$, there is an infinite number of solutions
$$                                                                                                      U(x,y)=\left\{ \begin{array}{ll} \displaystyle                                                                         \frac{KR}{2}(\frac{r^2}{a^2}-1),\quad\mbox{if}\ 0\leq r\leq a, \\ \\
\displaystyle KR\ln\frac{r}{a},\quad\mbox{if}\ a\leq r\leq R,
\end{array}\right.
$$
that change equation type for any $a$ ($0<a<R$). Implementing inequalities

$$
C<-\frac{KR}{2},\ \mbox{if}\ K\geq0;\ C<0,\ \mbox{if}\ K<0
$$
there is an infinite number of hyperbolic trivial solutions
\begin{equation}
U=\frac{Kr^2}{2R}+C.
\end{equation}
\textbf{The model Neumann problem for equation (2).} The problem has an infinite number of solutions
$$                                                                                                  U(x,y)=\left\{ \begin{array}{ll} \displaystyle                                                                         \frac{KR}{2a^2}(r^2-a^2)+KR\ln\frac{a}{R}+C,\quad\mbox{if}\ 0\leq r\leq a, \\ \\
\displaystyle KR\ln\frac{r}{R}+C,\quad\mbox{if}\ a\leq r\leq R,
\end{array}\right.
$$
that change equation type and an infinite number of trivial hyperbolic solutions
$$
U=\frac{Kr^2}{2R}+C
$$
for any $C,a,\ a\in(0,R).$

The above model Dirichlet and Neumann problems for equations (1) and (2), thanks to random parameter $a$,
have an infinite number of solutions that change the type of the equation.
When Dirichlet and Neumann conditions for equation (1) combined, the value of $a$ is uniquely determined.

\textbf{Model case of Cauchy problem.} $D$ is a circle of radius $R$, $U|_{r=R}=H,\ \frac{\partial U}{\partial r}|_{r=R}=K.$\\
\textbf{Model case of Cauchy problem for equation (1).}

Let $H>0,\ K>0.$ Satisfying Neumann condition in solution (5) of Dirichlet problem
 $\frac{\partial U}{\partial r}|_{r=R}=K,$
is determined value $a:$
$$
a=Re^{-\frac{H}{KR}}.
$$
We have one solution of the model Cauchy problem changing the type of the equation.

Satisfying in trivial hyperbolic solutions (8) of the model Neumann the Dirichlet condition,
we obtain a single trivial hyperbolic model solution of the Cauchy problem :
$$
U=\frac{K}{2R}(r^2-R^2)+H,\quad H\leq0, \ \mbox{èëè}\
H-\frac{KR}{2}<0,
$$
where $K$ of any sign.\\
\textbf{Model case of Cauchy problem for equation (2).} The problem has an infinite number of solutions
$$                                                                                                  U(x,y)=\left\{ \begin{array}{ll} \displaystyle                                                                         \frac{KR}{2a^2}(r^2-a^2)+KR\ln\frac{a}{R}+H,\quad\mbox{if}\ 0\leq r\leq a, \\ \\
\displaystyle KR\ln\frac{r}{R}+H,\quad\mbox{if}\ a\leq r\leq R,
\end{array}\right.
$$
that change equation type and an infinite number of trivial hyperbolic solutions
$$
U=\frac{Kr^2}{2R}+H-\frac{KR}{2}
$$
for any $a,\ a\in(0,R).$

Further only problems for equation (1) are considered.

\textbf{The model Dirichlet problem for the upper half-plane.}
Domain $D$ is the upper half plane $y>0$, $U|_{y=0}=0.$
The problem has an infinite number of solutions that change the type of equations:
$$
U(x,y)=\left\{ \begin{array}{ll}
\displaystyle Ce^x\sin y,\quad\mbox{if}\ 0\leq y\leq\pi,\\
\displaystyle Ce^x sh(\pi-y),\quad\mbox{if}\ y\geq\pi,
\end{array}\right.
$$
for any  $C>0$. When $C=0$ there is the zero solution.

\textbf{The solution on the whole plane.} The function
$$                                                                                                       U(x,y)=\left\{ \begin{array}{ll} \displaystyle                                                                         C(r^2-a^2),\quad\mbox{if}\ 0\leq r\leq a, \\ \\
\displaystyle 2Ca^2\ln\frac{r}{a},\quad\mbox{if}\ a\leq r<
+\infty,
\end{array}\right.
$$
for arbitrary $a>0,\ C>0$ is a solution of the equation (1) that changes the type of the equation on the plane.

\textbf{The Goursat problem}. In the domain $D:\ y+x> 0,\ y-x+1> 0$ it is required to find a solution to
equation (1) satisfying the conditions
\begin{equation}
U|_{y=-x}=f_1(x),\ x\leq\frac{1}{2};\quad U|_{y=x-1}=f_2(x),\
x\geq\frac{1}{2};\quad f_1(\frac{1}{2})=f_2(\frac{1}{2}).
 \end{equation}
Lines $y=-x,\ y=x-1$ are the characteristics of the wave equation. The function
\begin{equation}
U= f_1(\frac{x-y}{2})+f_2(\frac{x+y+1}{2})-f_1(\frac{1}{2})
\end{equation}
is the solution of the wave equation that satisfies (9).

If $(x,y)\subset D,$ then $\frac{x-y}{2}\leq\frac{1}{2},\
\frac{x+y+1}{2}\geq\frac{1}{2}.$ Requiring, for example,
$$
f_1(x)-f_1(\frac{1}{2})\leq0,\ f_2(x)\leq0 ,$$ we obtain that
(10) is negative in $D$, and thus it is a trivial solution of the hyperbolic Goursat problem.

\textbf{Tricomi problem}.  Let us consider the classical formulation of the Tricomi problem
for the Lavrentev-Bitsadze equation (3). Let $D$ be a domain bounded by the line $\sigma$ with endpoints
$A(-1,0),\ B(1,0),$ located in the upper $y>0$ and with the characteristics of
 $AC\ y=-x-1,BC\ y=x-1$  if $y<0$. $D$ is the characteristic domain of the Tricomi problem.
 It is required to determine the function $U$ with the following properties:

1) $U(x,y)$ is the solution of equation (1) in $D$ if
$y\neq0$;

2) is continuous in the enclosed domain $D$;

3) the partial derivatives $ \frac{\partial U}{\partial x},\
\frac{\partial U}{\partial y}$ are continuous within the region $D$,
and near the points $A,\ B$ , they can turn into infinity order less than one;

4) on the line $\sigma$ and on the characteristics $AC$ takes setpoints

\begin{equation}
 U|_\sigma=\varphi(s),
\end{equation}
\begin{equation}
 U|_{AC}=f(x).
\end{equation}

In A.V. Bitsadze [2] they proved the existence and uniqueness of the solution of the Tricomi problem reviewed above, and if the boundary of $\sigma$ is a semicircle, the solution is obtained explicitly.

Let $D_1,\ D_2$  be parts of the domain $D$, where $y>0$ and $y<0$, respectively.

We show that there exists a solution of the equation (1) that changes its type and satisfies only one condition
\begin{equation}
 U|_\sigma=\varphi(s)\geq0,\quad\ \varphi(-1,0)=\varphi(1,0)=0,
\end{equation}
which is positive in $D_1$, negetive in $D_2$ and does not require the boundary condition (12) on the characteristics of $AC$.

Let $U(x,y)$ be harmonic function in $D_1$ satisfying (13) and the condition
\begin{equation}
 U(x,0)=0,\quad -1\leq x\leq 1.
\end{equation}
 By the maximum principle and the boundary conditions (13), (14), the harmonic function
 $U(x,y)$ is positive in $D_1$. Since at
$y=0$ it has the lowest value, then, by the Hopf lemma for harmonic functions
 $\frac{\partial U}{\partial n}|_{y=0}<0$ for all $x\in(-1,1)$. The normal $n$ is external to $D_1$. Hence
$$
 \left.\nu(x)=\frac{\partial U}{\partial y}\right|_{y=0}>0,\quad -1<x<1.
 $$
 The function
\begin{equation}
 U(x,y)=\frac{1}{2}\int\limits_{x-y}^{x+y} \nu(x)dx
\end{equation}
is the solution of the wave equation in $D_2$ and
$$
\left.U(x,0)=0,\quad \frac{\partial U}{\partial
y}\right|_{y=0}=\nu(x),\quad -1<x<1.
$$
Since $\nu(x)>0$ and in $D_2$ $x+y<x-y$, then the solution given by (15) is negative in $D_2$.
Setting in (15)
$x+y=-1$, we find the value of the solution to the characteristic $AC$:
$$
f(x)=U|_{AC}=-\frac{1}{2}\int\limits_{-1}^{2x+1} \nu(x)dx<0,\
x\in(-1,0],\ U(-1,0)=0.
$$

Thus, setting only one boundary condition (13), we have a solution of the equation (1)
that changes its type in the domain, which is characteristic to the Tricomi problem.

Note that the search for the solution was begun from satisfying the boundary condition (13),
at first finding a solution in $D_1$ ($y>0$), and then in
$D_2$ ($y<0$). In other words, moving \textit{\textbf{"top - down".}}

If $\sigma$ is a semicircle: $x^2+y^2=1,\ y\geq0$, the solution can be written explicitly.
In $D_1$:
$$
U(x,y)=\frac{1}{\pi}\int\limits_{\sigma}\varphi\frac{\partial
G}{\partial n}ds,
$$
$$
G(\bar{x},\bar{\xi})=\ln\frac{|\bar{\xi}^*-\bar{x}||\bar{\xi}^--\bar{x}|}{|\bar{\xi}^+-\bar{x}||\bar{\xi}-\bar{x}|},\,
\bar{x}=(x,y),\,\bar{\xi}=(\xi,\eta),\,
\bar{\xi}^*=\frac{1}{|\bar{\xi}|}\bar{\xi},\,\bar{\xi}^-=(\xi,-\eta),\,
\bar{\xi}^+=\frac{1}{|\bar{\xi}|}\bar{\xi}^-.
$$
$G(\bar{x},\bar{\xi})$~--- Green's function of the Dirichlet problem for the domain $D_1.$

After direct calculation:
$$
\nu(x)=\frac{\partial U(x,0)}{\partial
y}=\frac{4(1-x^2)}{\pi}\int\limits_0^\pi \varphi(\cos,\sin\theta)
\frac{\sin\theta}{(1-2x\cos\theta+x^2)^2}d\theta.
$$
Solution of the problem in $D_2$ is defined by (15).

Let us find the solution of equation (1), which positive in  $D_1$ and negative in
$D_2$, but, when moving \textit{\textbf{"bottom - up"}},
satisfying the first boundary condition (12) on the characteristics of $AC$.

We seek the solution in $D_2$ in the form
$$
U(x,y)= f_1(x+y)+ f_2(x-y).
$$
Satisfying the condition (12) on the characteristics of $AC (y=-x-1)$ and the condition (14), we obtain
\begin{equation}
U(x,y)=f(\frac{x-y-1}{2})-f(\frac{x+y-1}{2}).
\end{equation}
In the domain $D_2:$ $\frac{x-y-1}{2}>\frac{x+y-1}{2}.$
Solutions in $D_2$ will be negative, if the function $f(x)$
is decreasing on $(-1,0)$. Taking into account that $f(-1)=0$,
we obtain the condition of negativity in $D_2$ solutions (16):
$$
f(-1)=0,\ \mbox{function}\ f(x)\ \mbox{decreasing,}\ f(x)<0,\ -1\leq
x\leq0.
$$

From (16):
\begin{equation}
\nu(x)=\frac{\partial U(x,0)}{\partial
y}=-f^\prime(\frac{x-1}{2}).
\end{equation}
For the existence of $\nu(x)=\frac{\partial U(x,0)}{\partial y}$
it is required the existence of $f^\prime(x)$. Since $f(x)$ is decreasing,
 then $f^\prime(x)<0$. Then, according to (17)
$$
\nu(x)=\frac{\partial U(x,0)}{\partial y}>0,\ -1<x<1.
$$

To find a solution in $D_1$ with the condition (14) we come to the problem of finding in $D_1$
a harmonic function satisfying
$$
U|_\sigma=\varphi(s)\geq0,\quad U(x,0)=0,\quad \frac{\partial
U(x,0)}{\partial y}=-f^\prime(\frac{x-1}{2})>0.
$$
The obtained problem is overdetermined, if the domain
$D_1$ or boundary function $\varphi(s)$, is pre-defined, not even requiring non-negativity of $\varphi(s)$.
If in a certain neighborhood of $y>0$
there exists a harmonic function, which is the solution of the Cauchy problem
\begin{equation}
 U(x,0)=0,\quad
\frac{\partial U(x,0)}{\partial y}=-f^\prime(\frac{x-1}{2}),\quad
-1<x<1,
\end{equation}
then the function $\nu(x)=-f^\prime(\frac{x-1}{2})$ is analytic on $x$ on the interval $(-1,1).$
It follows that the function given on the characteristic $AC$ must be analytic.

By Kowalewski theorem, in a neighborhood $y>0,\ -1<x<1$
there exists a harmonic function satisfying the analytical conditions of the Cauchy problem (18).
From the conditions
$$
 U(x,0)=0,\quad
\frac{\partial U(x,0)}{\partial
y}=-f^\prime(\frac{x-1}{2})>0,\quad -1<x<1
$$
follows the existence of a neighborhood of $y>0$, which is adjacent to the entire interval
$(-1,1)$, where this harmonic function is positive.
This neighborhood, or any other one contained in it,
which is also adjacent to the entire interval $(-1,1)$,
defines $D_1$ when formulating of the Tricomi problem,
if setting the boundary value only on the characteristics of $AC$:
\begin{equation}
U|_{AC}=f(x),\   f(-1)=0,\ f^\prime(x)<0,\ \mbox{function}\ f(x)\
\mbox{is analytical.}
\end{equation}

We obtain an explicit formula of the solution of the Cauchy problem (18) for the Laplace equation.
Turning to the complex plane $ z=x+\imath y,\quad
\overline{z}=x-\imath y, $ Laplace equation can be written as
$$
\frac{\partial^2 U}{\partial z\partial \overline{z}}=0.
$$
Taking
$$
U= f_1(z)+ f_2(\overline{z}),
$$
after satisfying the conditions (18) of the Cauchy problem, we obtain
\begin{equation}
U(x,y) =\imath(f(\frac{z-1}{2})-f(\frac{\overline{z}-1}{2}))
\end{equation}
The area of harmonicity is defined by the possibility of analytical extension of the function
$f(\frac{x-1}{2})$ from the interval $(-1,1)$.

Thus, we come to the two structures (area + boundary conditions) of, perhaps,
many formulations of boundary value problems for the equation (1).

Let $D$ be a domain that is characteristic to the Tricomi problem.\\
\textbf{Tricomi problem 1.} The boundary condition (13) is defined only on $\sigma$ the upper boundary of $D_1$.\\
\textbf{Tricomi problem 2.} Only one condition (12) is set on the characteristics of $AC$.
The domain $D_1$ is defined during the solution process.

Consider few examples of the Tricomi problem 2.

\textbf{Example.} $f(x)=-(1+x)^3$. In accordance with (16), (20)
$$
U(x,y)=\left\{ \begin{array}{ll} \displaystyle
\frac{y}{4}(3(x+1)^2-y^2),\quad\mbox{if}\ (x,y)\in D_1, \\ \\
\displaystyle \frac{y}{4}(3(x+1)^2+y^2),\quad\mbox{if}\ (x,y)\in
D_2.
\end{array}\right.
$$
$D_1$ is an arbitrary domain in the Tricomi problem, points of which satisfy
$\sqrt{3}(x+1)-y\geq0.$ In this case
$U(x,y)>0$ in $D_1$.

\textbf{Example.} Let the boundary function $f(x)$, defined on
$[-1,0]$, be represented by a series
$$
f(x)=\sum_{n=0}^\infty a_n(x+\frac{1}{2})^n.
$$
In accordance with (16), in $D_2$:
$$
U=\sum_{n=0}^\infty \frac{a_n}{2^n}((x-y)^n-(x+y)^n).
$$
In accordance with (20), in $D_1$:
\begin{equation}
U=\imath\sum_{n=0}^\infty
\frac{a_n}{2^n}(z^n-\overline{z}^n)=-\sum_{n=0}^\infty
\frac{a_n}{2^{n-1}}r^n\sin(n\varphi).
\end{equation}
Domain $D_1$ is determined by examining the sign of the function (21).

\textbf{Note 1}. It should be noted that in all the tasks
there remains an open question whether all solutions are found or not,
even if the problem got an infinite number of solutions.

\textbf{Note 2}. Equations (1), (3) as well as Laplace and wave equations are the special cases of equation (2).
Thus, all boundary problems typical to Laplace, wave, Lavrentev-Bitsadze equations and considered above for equation (1)
are included in possible statements of boundary problems for equation (2).

\textbf{Note 3.} Examples of formulations of boundary value problems Tricomi 1 and 2
give an opportunity to suppose that two components,
which define the boundary problem (domain and boundary conditions),
may be produced according to the method of finding a solution, that is,
setting the boundary value problem depends on the method of solving it.

\textbf{Note 4.}  Tricomi problem 2 shows a particular impact of the domain's geometry to the formulation of boundary value problems for the equation (1). Let $D$ be the domain domain that is characteristic of the Tricomi.
If, when formulating the Dirichlet problem for Lavrentev-Bitsadze equation, we come to the classical Tricomi problem, in which for one of the characteristics the boundary condition is not specified, then for the equation (1) with the additional condition of non-negativity of the function that is defined on the boundary of $\sigma$, the boundary conditions might not be already set on two characteristics.\\

\section*{2. Other equations.}
A natural extension of equation (2) is
\begin{equation}
\prod_{i=1}^n(L_i(u)-f_i(u))=0,
\end{equation}
where $L_i(u)$, for example, are linear differential operators of second order.
We look for a continuously differentiable in $D$
function and disjoint areas $D_i, \ \cap_{i=1}^n D_i=D,$ where
$L_i(u)=f_i(u).$ On the boundaries of $D$ and $D_i$ there specified the additional conditions.

If $L_i(u)=\Delta u$, then the equation (22) describes the general problem of matching
the vortex and potential flows of an ideal fluid [3].
Here $u(x,y)$ is a stream function, $\omega_i=f_i(u)$ is a vorticity in $D_i.$

 For example, if $D$ is an enclosed domain and a problem statement is
$$
\Delta U(x,y)(\Delta U(x,y)-\omega)=0,\ \omega>0,
U|_\Gamma=\varphi(s)\geq0,\ U|_{\Gamma_1}=0,
$$
$\Gamma_1$ is the boundary of $D_1\subset D,$ where $\Delta
U(x,y)=\omega,$ is the Goldshtik problem that describes separated flows by Lavrentev scheme [4].

The considered equations can be model in gas dynamics and hydromechanics problems,
heat,wave and other processes, which solutions have to satisfy different equations,
depending on processes settings and solution properties.

\markright{
 {\footnotesize   Isaac I.Vainshtein
 \hfill Model problems for two equations, which type depends on solution.}}




\bigskip




\begin{thebibliography}{4}

\bibitem{Vain}
 {\it Vainstein I.I.} The model problems for nonlinear analog Lavrent'ev-Bitsadze equations.
 International conference on the 80th anniversary of Academician Mikhail Lavrentiev.
 Inverse and ill-posed problems of mathematical physics. Abstracts. Novosibirsk, 2012, 352 p.
\bibitem{Vain}
 {\it Bitsadze A.V.} The equations of mixed type. Moscow, Academy of Science USSR Publishing house, 1959, 164 p.
\bibitem{Vain}
{\it Vainstein I.I.} Solution of two dual problems
of matching the vortex and potential flows by Goldshtik variational method.
Journal of Siberian Federal University. Ser. Mathematics and physics. volume 4, 2011, p. 320-331.
\bibitem{Vain}
 {\it Goldshtik M.A.} Vortical flows, Science, Novosibirsk, 1961, 364 p.
\end{thebibliography}
\end{document}